\documentclass{article}
\usepackage{graphicx}
\usepackage{amsmath}

\newtheorem{theorem}{Theorem}

\newtheorem{corollary}[theorem]{Corollary}

\newtheorem{lemma}[theorem]{Lemma}

\newtheorem{problem}[theorem]{Problem}

\newenvironment{proof}[1][Proof]{\textbf{#1.} }{\ \rule{0.5em}{0.5em}}

\begin{document}

\title{Isotopies of planar compacta need not extend to isotopies of planar continua}
\author{Paul Fabel \\
Mississippi State University}
\maketitle

\begin{abstract}
We construct, as a step toward a theory of braids over planar compacta, an
isotopy of a planar compactum which cannot be extended to an isotopy of any
planar continuum.
\end{abstract}

\section{\protect\bigskip Introduction}

In this paper we address an obstruction to developing a satisfactory theory
of braids over planar compacta by establishing a technical sense in which
all braids over planar continua are not wild. We then construct an example
that is apparently the first of its kind: An isotopy of a planar compactum
shown not to be the restriction of an isotopy of any planar continuum.

Given a finite set $F\subset R^{2}$ with $\left| F\right| =n$, Artin's pure
braid group on $n$-strands $G_{n}$ is $\pi _{1}(I(F,R^{2}),id_{F}),$ where $%
I(F,R^{2})$ is the space of embeddings of $F$ into $R^{2}.$ Artin's braid
group on $n$ strands is $\pi _{1}(Y,F)$ where $Y$ is the space of planar
sets of size $n$ in the Hausdorf metric. The simplest cases where $\left|
F\right| =\infty $ are still being explored. In \cite{Fabel1} the author
constructs a group of wild braids $\overline{B_{\infty }}$ which completes $%
B_{\infty },$ Artin's braid group on infinitely many strands. The full
inverse limit of Artin's pure braid groups is a group $G\subset \overline{%
B_{\infty }}$ consisting of the pure braids on infinitely many strands.
However, a typical braid representative $g\in \lbrack g]\in G$ is not the
restriction of an ambient isotopy of the plane.

In hopes of eliminating this particular brand of wildness, suppose instead $%
F\subset R^{2}$ is the disjoint union of finitely many nonseparating planar
continua $X_{1}\cup ..X_{n}$, and a braid representative is an isotopy $%
h_{t}:$ $F\hookrightarrow R^{2}$ such that $h_{0}=id_{F}$ and $h_{1}(F)=F.$
Though indirectly settled in case $F$ is PL (see \cite{Yag1}), a thorough
understanding of this braid group hinges on the following question:

\begin{problem}
\label{prob}Can every isotopy starting at $id_{W}$ of a planar continuum $W$
be extended to an isotopy of the plane?
\end{problem}

\bigskip We do not settle this challenging problem but prove in section \ref
{mainsect} Theorem \ref{nostretch}, a weaker technical result which says
roughly, that loops in the complements of finite subsets of the continuum $W$
cannot be forcibly stretched by large amounts via isotopy of $W.$ We then
construct in Corollary \ref{maincor} a wild isotopy over a convergent
sequence which is unextendable to an isotopy of any planar continuum.

The proof of Theorem \ref{nostretch} relies in part on Theorem \ref{lifting}%
, a result in section \ref{liftsect} which establishes mild conditions under
which homotopy classes in $X\backslash A$ are naturally transported via
homotopy of $A.$ In particular Theorem \ref{lifting} eliminates reliance on
canonical parameterizations of $X\backslash A.$

The final section \ref{lemmsect} includes Lemmas useful for the proof of
Theorem \ref{nostretch}.

\section{\label{liftsect}Theorem \ref{lifting}: a path lifting theorem.}

\subsection{Motivation}

Given an isotopy $h_{t}:A\hookrightarrow R^{2}$ and $[f]\in \pi
_{1}(R^{n}\backslash A),$ we would like $[f]$ to `come along for the ride'
via a natural choice $[f_{t}]\in \pi _{1}(R^{n}\backslash A_{t}).$ The
choice is clear in the special case $A$ is finite since we can canonically
extend $h_{t}$ to an isotopy of the plane. However we would like another
scheme for selecting $[f_{t}]$ since it's generally impossible to
canonically parameterize $X\backslash A_{t}$ in an arbitrary space $X.$
Fortunately we can generally recover $[f_{t}]$ under weak hypothesis as seen
in Theorem \ref{lifting}.

\subsection{$\mathcal{B}$ and $\mathcal{E}$}

\bigskip Suppose $X$ is a metric space and $\mathcal{B}$ is a collection of
compact subsets of $X$ such that $\forall A\in \mathcal{B}$ there exists a
compact $F\subset X$ and a homotopy $H:X\backslash A\times \lbrack
0,1]\rightarrow X\backslash A$ such that $\forall x\in X\backslash A$ $%
H(x,0)=x$ and $H(x,1)\in F.$

Decide $A_{n}\rightarrow A\in \mathcal{B}$ if $A_{n}\rightarrow A$ in the
Hausdorf metric and there exists $N$ and homeomorphisms $h_{n}:X\backslash
A\rightarrow X\backslash A_{n}$ for $n\geq N$ such that $h_{n}\rightarrow id$
in the compact open topology. For $\{A,B\}\subset \mathcal{B}$ let $\mathcal{%
H}(A,B)$ denote the Hausdorf distance.

Fixing a compact metric space $Z$ let $\mathcal{P}(A,f)$ denote the path
component of $f$ in $C(Z,X\backslash A).$ Let $\mathcal{E}$ denote the space
of ordered pairs $(A,\mathcal{P}(A,f))$ such that $A\in \mathcal{B}$ and $%
f\in C(Z,X\backslash A).$ Decide that the sequence $(A_{n},\mathcal{P}%
(A_{n},f_{n}))\rightarrow (A,\mathcal{P}(A,f))$ iff $A_{n}\rightarrow A$ in $%
\mathcal{B}$ and if there exists $g\in \mathcal{P}(A,f)$ and $N$ such that $%
g\in \mathcal{P}(A_{n},f_{n})$ $\forall n\geq N.$ Define $\Pi _{1}:\mathcal{%
E\rightarrow B}$ as $\Pi _{1}(A,\mathcal{P}(A,f))=A.$ Let $\Pi _{2}(A,%
\mathcal{P}(A,f))=\mathcal{P}(A,f).$ Topologize $\mathcal{B}$ and $\mathcal{E%
}$ by sequential convergence declaring a set $\mathcal{F}$ to be closed if
and only the limit of each convergent sequence in $\mathcal{F}$ is an
element of $\mathcal{F}$.

\begin{lemma}
\label{ezman}Suppose $\alpha :[0,1]\rightarrow \mathcal{B}$ is continuous $%
T\in \lbrack 0,1]$, $\varepsilon >0,$ and $F\subset X\backslash \alpha (T)$
is compact. Then there exists $\gamma >0$ so that if $\left| t-T\right|
<\gamma $ then $\mathcal{H}(\alpha (t),\alpha (T))<\varepsilon $ and there
exists a homeomorphism $h_{t}:X\backslash \alpha (T)\rightarrow X\backslash
\alpha (t)$ such that $d(h_{t}(x),x)<\varepsilon $ $\forall x\in F.$
\end{lemma}

\begin{proof}
Suppose $t_{n}\rightarrow T.$ Then $\mathcal{H}(\alpha (t_{n}),\alpha
(T))\rightarrow 0.$ Moreover there exists $N$ and homeomorphisms $%
h_{n}:X\backslash \alpha (T)\rightarrow X\backslash \alpha (t_{n})$ such
that $h_{n|F}\rightarrow h_{|F}$ uniformly for $n\geq N.$ If the conclusion
of the Lemma were false we could choose $t_{n}\rightarrow T$ so as to
contradict one of the previous two sentences.
\end{proof}

\begin{lemma}
\label{ezman2}Suppose $\beta :[0,1]\rightarrow \mathcal{E}$ is continuous
and $g\in \Pi _{2}(\beta (T)).$ Then there exists $\gamma >0$ such that if $%
\left| T-t\right| <\delta $ then $g\in \Pi _{2}(\beta (t)).$
\end{lemma}

\begin{proof}
Let $T\in \lbrack 0,1].$ Suppose in order to obtain a contradiction that
there does not exist $\gamma >0$ so that $g\in \Pi _{2}\beta (t)$ if $\left|
t-T\right| <\gamma .$ Choose $t_{n}\rightarrow T$ with $g\notin \Pi
_{2}\beta (t_{n})\forall n.$ Since $\beta (t_{n})\rightarrow \beta (T)$
choose $M$ and $f_{T}\in \Pi _{2}(\beta (T))$ so that if $n\geq M$ then $%
f_{T}\in \Pi _{2}(\beta (t_{n})).$ Let $\beta (t)=(A_{t},\mathcal{P}%
(A_{t},f_{t})).$ Let $H:Z\times \lbrack 0,1]\rightarrow X\backslash A_{T}$
be any homotopy between $g$ and $f_{T}.$ Since $im(H)$ is compact choose $%
\varepsilon >0$ so that each point of $im(H)$ is at least $\varepsilon $
from each point of $A_{T}$. $.$Choose by Lemma \ref{ezman} $N\geq M$ so that
if $n\geq N$ then $\mathcal{H}(A_{T},A_{t_{n}})<\varepsilon .$ Thus $H$ is a
homotopy between $f_{T}$ and $g$ in $X\backslash A_{t_{n}}.$ Hence $g\in \Pi
_{2}(\beta (t_{n}))$ for $n\geq N$ and we have a contradiction.
\end{proof}

\begin{theorem}
\label{lifting}Suppose $\alpha :[0,1]\rightarrow \mathcal{B}$ is continuous, 
$e\in \mathcal{E}$ and $\Pi _{1}(e)=\alpha (0).$ Then there exists a unique
continuous function $\beta :[0,1]\rightarrow \mathcal{E}$ such that $\alpha
=\Pi _{1}(\beta ).$
\end{theorem}

\begin{proof}
Let $I=\{t\in \lbrack 0,1]$ such that there exists a continuous function $%
\beta :[0,t]\rightarrow \mathcal{E}$ such that $\alpha _{\lbrack 0,t]}=\Pi
_{1}(\beta )\}$ and $\beta (0)=e\}.$ Note $I\neq \emptyset $ since $0\in I$
and $I$ is connected since $I$ is the union of intervals each of which
contains $0.$ To prove $I$ is open in $[0,1]$ suppose $T\in I$ and $T\in
\lbrack 0,1).$ Suppose $\beta :[0,T]\rightarrow \mathcal{E}$ satisfies $%
\alpha _{\lbrack 0,T]}=\Pi _{1}\beta .$ Let $f\in \Pi _{2}\mathcal{\beta }%
(T).$ Since $\alpha (T)$ and $im(f)$ are disjoint compacta choose $%
\varepsilon >0$ so that each point of $\alpha (T)$ is at least a distance $%
\varepsilon $ from each point of $im(f).$ Choose by Lemma \ref{ezman} $%
\delta >0$ so that if $\left| t-T\right| <\delta $ then $\mathcal{H}(\alpha
(t),\alpha (T))<\varepsilon /2.$ Hence $im(f)\subset X\backslash \alpha (t)$
and $\beta $ can be continuously extended to $[0,T+\delta )$ by defining $%
\beta (t)=(\alpha (t),\mathcal{P}(\alpha (t),f))$ for $t\in \lbrack
T,T+\delta ).$

To see that $I$ is closed suppose $[0,T)\subset I.$ Choose a compact $%
F\subset X\backslash \alpha (T)$, a homotopy $H:X\backslash \alpha (T)\times
\lbrack 0,1]\rightarrow X\backslash \alpha (T)$ such that $\forall x\in
X\backslash A$ $H(x,0)=x$ and $H(x,1)\in F.$ Since $\alpha (T)$ and $F$ are
disjoint compacta choose $\varepsilon >0$ so that each point of $\alpha (T)$
is at least a distance $\varepsilon $ from each point of $F.$ Choose by
Lemma \ref{ezman} $S<T$ and a homeomorphism $h_{S}:X\backslash \alpha
(T)\rightarrow X\backslash \alpha (S)$ so that if $S\leq t\leq T$ then $%
\mathcal{H}(\alpha (t),\alpha (T))<\varepsilon /2$ and $d(h_{S}(x),x)<%
\varepsilon /2$ $\forall x\in F.$ Let $\beta :[0,S]\rightarrow \mathcal{E}$
satisfy $\Pi _{1}\beta =\alpha _{\lbrack 0,S]}.$ Let $f\in \Pi _{2}\beta
(S). $ Define for $s\in \lbrack 0,1]$ $g^{s}\in C(Z,X\backslash \alpha (S))$
via $g^{s}(z)=h_{S}H(h_{S}^{-1}f(z),s).$ Note $g^{0}=f$ and hence $g^{1}\in 
\mathcal{P}(\alpha (S),f)=\Pi _{2}(\beta (S)).$ Note $im(g^{1})\subset
h_{S}(F).$ Hence $im(g^{1})\cap \alpha (t)=\emptyset $ for $t\in \lbrack
S,T].$ Defining $\gamma (t)=(\alpha (t),\mathcal{P}(\alpha (t),g^{1}))$ for $%
t\in \lbrack S,T]$ creates a continuous lift $\gamma \cup \beta _{\lbrack
0,S]}$ of $\alpha _{|[0,T]}$ and proves $T\in I.$ Hence $I=[0,1]$ and the
existence of $\beta $ is established.

For uniqueness suppose $\gamma $ and $\beta $ are two lifts of $\alpha $
such that $\beta (0)=\gamma (0)=e.$ Let $J=\{j\in \lbrack 0,1]|\gamma
(j)=\beta (j)\}.$ Note $J$ is nonempty and closed. It suffices to prove $J$
is open. Suppose $T\in J.$ Since $\gamma (T)=\beta (T)$ let $g\in \Pi
_{2}(\beta (T))\cap \Pi _{2}(\gamma (t)).$Choosing $\gamma $ as in Lemma \ref
{ezman2} it follows that $\beta (t)=\gamma (t)=(\alpha (t),\mathcal{P}%
(\alpha (t),g))$ for $\left| t-T\right| <\gamma .$ Hence $J$ is open and $%
\gamma =\beta .$
\end{proof}

\section{\label{mainsect}Main results: Theorem \ref{nostretch} and Corollary 
\ref{maincor}}

Let $\mathcal{B}$ denote the space of planar compacta, topologized as in
section \ref{liftsect}, each of which is the disjoint union of finitely many
nonseparating planar continua. ( i.e. $X\in \mathcal{B}$ iff $\exists n$ and
nonseparating planar continua $X_{1},..X_{n}$ such that $X=\cup
_{i=1}^{n}X_{i}$ and $X_{i}\cap X_{j}=\emptyset $ whenever $i\neq j$).

Fixing $Z=S^{1},$ define $\mathcal{E}$ as in section \ref{liftsect}$.$ For $%
\delta >0$ let $B(p,\delta )\subset R^{2}$ and $int(B(p,\delta ))\subset
R^{2}$ denote respectively the closed and open disks centered at $p$ of
radius $\delta .$ Let $\partial B(p,\delta )$ denote $B(p,\delta )\backslash
intB(p,\delta ),$ the circle of radius $\delta $ centered at $p.$

Given $A\subset R^{2},$ an \textbf{isotopy }$h_{t}:A\hookrightarrow R^{2}$
is a homotopy such that $h_{0}=id_{A}$ and $h_{t}$ is an embedding. We do
not require that $h_{t}(A)=A.$ Given an isotopy of a planar continuum $W,$
Theorem \ref{nostretch} establishes a technical sense in which small loops
in the complements of finite subsets of $W$ cannot be forcibly stretched by
large amounts.

Care is required since in general a small isotopy $h_{t}:A\hookrightarrow
R^{2}$ of a finite set $A\subset R^{2}$ can dramatically stretch
complementary loops when $h_{t}$ is extended to the plane. However if the
loop $S\subset R^{2}\backslash A$ encloses only one point of the finite set $%
A$, then an isotopy of $A$ cannot not stretch $S$ very much. The latter fact
is exploited in our proof of Lemma \ref{biglemma} very roughly sketched as
follows:

Assuming $W$ is nonseparating and the isotopy $h_{t}:W\hookrightarrow R^{2}$
fixes $p\in W$ throughout consider $X=B(p,\varepsilon _{1})\cup
(B(p,\varepsilon _{5})\cap W)$ with $Y_{1}\subset X$ the component
containing $p.$ Take a simple closed curve $S\subset R^{2}\backslash X$ such
that $Y_{1}$ is inside $S$ and the remaining components $Y_{2},...Y_{n}$of $%
X $ lie outside $S.$ Let $p=y_{1}$ and for $i\geq 2$ pick $y_{i}\in Y_{i}$
such that $\left| y_{i}-p\right| =\varepsilon _{5}$. Let $%
A=\{y_{1},y_{2},...y_{n}\}.$ Let $H_{t}:R^{2}\backslash X_{t}\rightarrow
R^{2}\backslash A_{t}$ be an $\varepsilon _{8}$ homeomorphism mapping a
neighborhood of $Y_{i}^{t}$ to a neighborhood of $y_{i}^{t}.$ Now suppose $%
F\subset W$ is finite and $f:S^{1}\rightarrow R^{2}\backslash F$ satisfies $%
im(f)\subset B(p,\varepsilon _{1}).$ Observe $S_{t}=H_{t}^{-1}(\partial
B(p,\varepsilon _{3}))$ has small diameter and encloses $im(f_{t})$. Hence $%
im(f_{t})$ has small diameter.

\begin{lemma}
\label{biglemma} Suppose $W\subset R^{2}$ is a continuum, $%
h_{t}:W\hookrightarrow R^{2}$ is an isotopy, $p\in W$, $h_{t}(p)=p\forall
t\in \lbrack 0,1]$ and $\varepsilon >0.$ Suppose $\forall n\geq 1$ $%
F_{n}\subset W$ is finite, $f_{n}:S^{1}\hookrightarrow R^{2}\backslash F_{n}$
is an embedding, $im(f_{n})\rightarrow \{p\}$ in the Hausdorf metric, $%
\alpha _{n}:[0,1]\rightarrow \mathcal{B}$ is defined as $\alpha
_{n}(t)=h_{t}(F_{n})$ and $\beta _{n}:[0,1]\rightarrow \mathcal{E}$ is the
(unique) lift of $\alpha _{n}$ such that $\Pi _{2}\beta _{n}(0)=\mathcal{P}%
(F_{n},f_{n}).$ Then there exists $N$ such that if $n\geq N$ and $T\in
\lbrack 0,1]$ then there exists $f_{n}^{T}\in \Pi _{2}\mathcal{\beta }%
_{n}(T) $ such that diam($im(f_{n}^{T})$)$<\varepsilon .$
\end{lemma}

\begin{proof}
Note first that $\beta _{n}$ is well defined by Theorem \ref{lifting}. By
uniform continuity of the isotopy $h$ over $[0,1]\times $ $W$, and by
injectivity of $h_{t}$, choose $0<\varepsilon _{1}<\varepsilon
_{2}....<\varepsilon _{10}=\varepsilon $ such that $\forall \{s,t\}\subset
\lbrack 0,1],$ $\forall \{x,y\}\subset W$ if $|h_{s}(x)-h_{s}(y)|\leq
\varepsilon _{n}$ then $\left| h_{t}(x)-h_{t}(y)\right| <\varepsilon _{n+1}$
and such that $\varepsilon _{n}<\varepsilon _{n+1}/2.$ Note if $W\subset
B(p,\varepsilon _{5})$ the theorem is proved since $h_{t}(W)\subset
int(B(p,\varepsilon _{6}))\subset B(p,\varepsilon )$ and we can canonically
deform all loops into $B(p,\varepsilon _{6})$ via radial contraction outside 
$R^{2}\backslash int(B(p,\varepsilon _{6}))$. So we assume $W\backslash
B(p,\varepsilon _{5})\neq \emptyset .$ Choose $N$ so that if $n\geq N$ then
im$(f_{n})$)$\subset int(B(p,\varepsilon _{1})).$ Suppose $n\geq N$. Let $%
X=B(p,\varepsilon _{1})\cup F_{n}\cup (B(p,\varepsilon _{5})\cap W).$ Let $%
P=X\cup V$ where $V$ is the union of the bounded complementary domains of $%
X. $ For clarity note

\begin{enumerate}
\item  Each component $P_{i}$ of $P$ is a nonseparating planar continuum of
diameter no greater than $2\varepsilon _{5}.$

\item  (By Lemma \ref{bd}) each component of $P_{i}$ of $P$ contains a point 
$y_{i}$ such that $\left| y_{i}-p\right| \geq \varepsilon _{5}.$
\end{enumerate}

Let $Y\subset P$ consist of those components $Y_{1},...Y_{K}$ of $P$ such
that $Y_{i}\cap F_{n}\neq \emptyset $ indexed so that $p\in Y_{1}.$ Since $%
2\varepsilon _{5}<\varepsilon _{6}$ each component of $h_{t}(Y)$ has
diameter no more than $\varepsilon _{7}.$ Hence we may choose $\forall t\in
\lbrack 0,1]$ a collection of disjoint closed topological disks $%
D_{1}^{t},...D_{K}^{t}$ such that diam$(D_{i}^{t})<\varepsilon _{8}$ and $%
int(D_{i}^{t})$ contains exactly one component of $h_{t}(Y),$ namely $%
h_{t}(Y_{i}).$ Choose $y_{i}\in Y_{i}$ such that $\left| p-y_{i}\right| \geq
\varepsilon _{5}$ for $i\geq 2$ $.$ Let $y_{1}=p.$ Let $A=%
\{y_{1},y_{2},...y_{K}\}.$ Choose $\forall t\in \lbrack 0,1]$ a surjective
`decomposition' map $H_{t}:R^{2}\rightarrow R^{2}$ such that, preserving
orientation, $H_{t}$ maps $R^{2}\backslash h_{t}(Y)$ homeomorphically onto $%
R^{2}\backslash h_{t}(A),$ $H_{t}$ maps $h_{t}(Y_{i})$ onto $%
\{h_{t}(y_{i})\} $ and $H_{t}$ fixes $R^{2}\backslash (\cup
_{i=1}^{K}int(D_{i}^{t}))$ pointwise. Let $r:S^{1}\rightarrow \partial
(B(p,\varepsilon _{3}))$ be any homeomorphism. Let $g_{n}^{t}=H_{t}^{-1}(r).$
Note $g_{n}^{t}$ is well defined since $h_{t}(A)\cap \partial
B(p,\varepsilon _{3})=\emptyset $. Note $g_{n}^{0}\in \mathcal{P}(A,r)$
since $A\cap B(p,\varepsilon _{3})=\{p\}.$ Hence $Y_{1}$ is contained within
the bounded complementary domain of $im(g_{n}^{0}).$ Moreover since $%
B(p,\varepsilon _{1})\subset Y_{1}$ it follows that im$(f_{n})$ is contained
within the bounded complementary domain of $im(g_{n}^{0}).$ By Lemma \ref
{conty} the map $\beta _{n}^{2}:[0,1]\rightarrow \mathcal{E}$ defined as $%
\beta _{n}^{2}(t)=(h_{t}(Y),\mathcal{P}(h_{t}(Y),g_{n}^{t}))$ is continuous.
Hence there exists $f_{n}^{T}\in \Pi _{2}\beta _{n}(T)$ such that im$%
(f_{n}^{T})$ is contained inside the bounded complementary domain of $%
im(g_{n}^{T}).$ In particular diam($f_{n}^{T}$)$\leq $diam($g_{n}^{T}$)$%
<2\varepsilon _{8}<\varepsilon _{9}<\varepsilon .$
\end{proof}

\begin{theorem}
\label{nostretch}Suppose $h_{t}:W\hookrightarrow R^{2}$ is an isotopy of the
continuum $W\subset R^{2}.$ Then for each $\varepsilon >0$ there exists $%
\delta >0$ such that if $F\subset W$ is finite, $f:S^{1}\hookrightarrow
R^{2}\backslash F$ is an embedding $\alpha :[0,1]\rightarrow \mathcal{B}$ is
defined as $\alpha (t)=h_{t}(F)$ and $\beta :[0,1]\rightarrow \mathcal{E}$
is the (unique) lift of $\alpha $ such that $\beta (0)=(F,\mathcal{P}(F,f))$
then $\forall T\in \lbrack 0,1]$ then there exists $f^{T}\in \Pi _{2}%
\mathcal{\beta }(T)$ such that diam($im(f^{T})$)$<\varepsilon .$
\end{theorem}

\begin{proof}
We reduce the problem to Lemma \ref{biglemma} as follows. If Theorem \ref
{nostretch} were false we could construct an isotopy $h_{t}:W\hookrightarrow
R^{2},$ $\varepsilon >0,$ sequences $t_{n}\in \lbrack 0,1],$ finite sets $%
F_{n}\subset W,$ embeddings $f_{n}:S^{1}\rightarrow R^{2}\backslash F_{n}$
such that diam$(imf_{n})\rightarrow 0$ and so that each element of $\Pi
_{2}(\beta _{n}(t_{n}))$ has diameter greater than $\varepsilon ,$ where $%
\beta _{n}$ is the unique lift guaranteed by Theorem \ref{lifting} of $%
\alpha _{n}:[0,1]\rightarrow \mathcal{B}$ defined as $\alpha
_{n}(t)=h_{t}(F_{n}).$ In such a counterexample each bounded complementary
domain of $im(f_{n})$ must include at least two points of $F_{n}\subset W$
to ensure each element of $\Pi _{2}(\beta _{n}(t_{n}))$ has diameter at
least $\varepsilon $. Since $W$ is compact, passing to a subsequence if
necessary, we may further assume $im(f_{n})\rightarrow p$ in the Hausdorf
metric. Finally if $h_{t}$ is a counterexample to the theorem then the
isotopy $g_{t}:W\hookrightarrow R^{2}$ defined as $%
g_{t}(w)=h_{t}(w)+p-h_{t}(p)$ is also a counterexample. Hence we may further
assume that $h_{t}(p)=p\forall t\in \lbrack 0,1].$ Now apply Lemma \ref
{biglemma} to conclude that no such counterexample exists.
\end{proof}

\begin{corollary}
\label{maincor}There exists an isotopy $h_{t}$ of a compactum $W\subset
R^{2} $ such that $h_{t}$ cannot be extended to an isotopy of any planar
continuum.
\end{corollary}

\begin{proof}
Consider $W=\{0\}\cup \{1,1/2,1/3...\}\subset R^{2}$ and the following
isotopy $h_{t}:A\hookrightarrow R^{2}.$ For $n\geq 2$ let the point $1/n$
perform a counterclockwise orbit around the point$\frac{1}{n-1}$ for $t\in
\lbrack 1/n,1/n-1]$ while all other points remain fixed. Taking $%
F_{n}=\{1/n,1/n-1,...1/2,1\}$, Let $f_{n}:S^{1}\hookrightarrow
R^{2}\backslash F$ be an embedding onto the boundary of a small convex set
containing $\{\frac{1}{n},\frac{1}{n-1}\}.$ Define $\alpha :[0,1]\rightarrow 
\mathcal{B}$ as $\alpha (t)=h_{t}(F_{n}).$ Now observe that each point of $%
\Pi _{2}\mathcal{\beta }(1)$ has diameter at least $1-\frac{1}{n-1}$ where $%
\beta :[0,1]\rightarrow \mathcal{E}$ is the unique lift of $\alpha $ such
that $\beta (0)=(F_{n},\mathcal{P}(F_{n},f)).$ Hence, for any continuum $W$
such that $A\subset W,$ we cannot extend $h_{t}$ to an isotopy of $W$, since
fixing $\varepsilon <1$ there does not exist $N$ guaranteed by Theorem \ref
{nostretch}.
\end{proof}

\section{\label{lemmsect}Lemmas}

\begin{lemma}
\label{bd}Suppose $X$ is a metric continuum, $p\in X$ and $\overline{B(p,1)}%
\neq X$ then each component of $\overline{B(p,1)}\cap X$ contains a point
whose distance from $p$ is exactly $1.$ ( $B(p,1)$ denotes the closed metric
ball of radius $1$ centered at $p.$)
\end{lemma}

\begin{proof}
To obtain a contradiction let $Y$ be a component of $B(p,1)\cap X$ violating
the conclusion. Since $X$ is compact, the components of $B(p,1)\cap X$ are
exactly the quasicomponents. Let $U$ and $V$ be nonempty open sets
separating $B(p,1)\cap X$ such that $Y\subset U.$ Now map the connected
space $X$ onto the disconnected space $U\cup \{\infty \}$ via $f(x)=x$ if $%
x\in U$ and $f(x)=\infty $ otherwise, and we have a contradiction.
\end{proof}

\begin{lemma}
\label{Alexander}Suppose $Y=\cup _{i=1}^{N}Y_{i}$ is the disjoint union of
finitely many nonseparating planar continua $Y_{1},..Y_{N}$ and $A\subset Y$
contains exactly one point from each continuum $Y_{i}.$ Then $\mathcal{F}%
_{Y} $ is nonempty and path connected where, in the compact open topology, $%
\mathcal{F}_{Y}$ is the space of maps $h:R^{2}\rightarrow R^{2}$ such that $%
h $ is an orientation preserving homeomorphism between $R^{2}\backslash Y$
and $R^{2}\backslash A$, $h$ maps $Y$ onto $A$, and there exists a sequence
of disjoint closed topological disks $D_{1},D_{2},...D_{N}$ ( which may vary
with $h$) such that $Y_{i}\subset int(D_{i})$ and $h$ fixes $A\cup
(R^{2}\backslash (\cup _{i=1}^{n}int(D_{N})))$ pointwise.
\end{lemma}

\begin{proof}
Fixing $f\in \mathcal{F}_{Y}$ the map $H:\mathcal{F}_{Y}\rightarrow \mathcal{%
F}_{A}$ defined as $H(h)=hf^{-1}$is a homeomorphism. Thus it suffices to
show $\mathcal{F}_{A}$ is path connected. Take any $g\in \mathcal{F}$ and
disjoint closed disks $D_{i}$ outside of which $g$ is pointwise fixed and
isotop $g$ to $id$ via the `Alexander isotopy' on each $D_{i}.$
\end{proof}

\begin{lemma}
\label{conty}Suppose $\alpha :[0,1]\rightarrow \mathcal{B}$ is continuous
and $\mathcal{B}$ is the collection of planar compacta $Y$ such that $Y$ is
the disjoint union of finitely many nonseparating planar continua $Y_{1}\cup
...\cup Y_{n}$. Suppose $\gamma :[0,1]\rightarrow \mathcal{B}$ is continuous
and $\gamma (t)$ contains exactly one point from each component of $\alpha
(t).$ Suppose $\forall t\in \lbrack 0,1]$ $H_{t}:R^{2}\backslash \alpha
(t)\rightarrow R^{2}\backslash \gamma (t))$ is a homeomorphism such that $%
H_{t}$ fixes pointwise an open set $U_{t}$ such that $R^{2}\backslash U_{t}$
is a collection of disjoint closed topological disks $D_{1}^{t},...D_{n}^{t}$
such that $int(D_{i}^{t})$ contains exactly one component of $\alpha (t).$
Suppose $Z=S^{1}$ and $r:S^{1}\rightarrow R^{2}$ is continuous and $\forall
t\in \lbrack 0,1]$ $im(r)\subset R^{2}\backslash \gamma (t).$ Then the map $%
\beta :[0,1]\rightarrow \mathcal{E}$ defined as $\beta (t)=(\alpha (t),%
\mathcal{P}(\alpha (t),H_{t}^{-1}r))$ is continuous.
\end{lemma}

\begin{proof}
Suppose $T\in \lbrack 0,1].$ By Lemma \ref{Alexander} $\beta (T)$ does not
depend on our choice of $D_{1}^{T},...D_{n}^{T}$ as long as $int(D_{i}^{T})$
contains exactly one component of $\alpha (T).$ Thus for $t$ sufficiently
close to $T$ $\beta (t)=(\alpha (t),\mathcal{P}(\alpha (t),H_{T}^{-1}r))$
and in particular $H_{T}^{-1}r\in \Pi _{2}\beta (t).$ Thus $\beta $ is
continuous.
\end{proof}

\end{document}